\documentclass[letterpaper, 10 pt, conference]{ieeeconf}  

\IEEEoverridecommandlockouts                              
\usepackage{graphics} 
\usepackage{amsmath} 
\usepackage{amssymb}  

\usepackage{xcolor}
\usepackage{tikz}
\usepackage{algorithm}
\usepackage{caption}  
\usepackage{booktabs}  
\usepackage{hyperref}
\usepackage[capitalise]{cleveref}  
\usepackage{url}  

\graphicspath{{code_figs/}}

\newtheorem{remark}{Remark}

\newtheorem{proposition}{Proposition}

\newcommand{\cC}{\mathcal{C}}
\newcommand{\cB}{\mathcal{B}}

\newcommand{\cCH}{\cC_\mathcal{H}}

\newcommand{\cCV}{\cC_\mathcal{V}}
\newcommand{\psiV}{\Psi_V}
\newcommand{\innerproduct}[2]{#1^T #2}

\newcommand{\norm}[2]{\|#1\|_{#2}}

\newcommand{\posvector}{\ensuremath{p}}

\newcommand{\rate}{\eta}
\newcommand{\numrays}{m}

\DeclareMathOperator{\interior}{int}

\DeclareMathOperator{\convexhull}{ch}
\DeclareMathOperator{\offdiag}{off-diag}

\DeclareMathOperator{\vect}{vec}

\newcommand{\tconeint}[1]{\interior \mathcal{K}_{#1}}

\title{\LARGE \bf
Polyhedral Lyapunov Functions with Fixed Complexity 
}

\author{Dimitris Kousoulidis and Fulvio Forni%
\thanks{D. Kousoulidis is supported by the Engineering and Physical Sciences Research Council (EPSRC) of the United Kingdom.}%
\thanks{D. Kousoulidis and F. Forni are with the Department of Engineering, University of Cambridge, CB2 1PZ, UK {\tt\small dk483@eng.cam.ac.uk \& f.forni@eng.cam.ac.uk}}
}

\begin{document}

\maketitle
\thispagestyle{empty}
\pagestyle{empty}

\begin{abstract}
Polyhedral Lyapunov functions can approximate any norm
arbitrarily well.
Because of this, they are used to study the stability
of linear time varying and linear parameter varying systems 
without being conservative.
However, the computational cost associated with using them
grows unbounded as the size of their representation increases.
Finding them is also a hard computational problem.

Here we present an algorithm that
attempts to find polyhedral functions while keeping
the size of the representation fixed,
to limit computational costs.  
We do this by measuring the gap from contraction
for a given polyhedral set.
The solution is then used to find perturbations
on the polyhedral set 
that reduce the contraction gap.
The process is repeated until a valid
polyhedral Lyapunov function is obtained.

The approach is rooted in linear programming.
This leads to a flexible method capable of handling
additional linear constraints and objectives,
and enables the use of the algorithm for control synthesis.
\end{abstract}

\section{Introduction}

In this paper we provide a tool for
the analysis and design of Linear Time Varying (LTV)
and Linear Parameter Varying (LPV) systems. 
We first focus on verifying the stability of these systems.
Then, we focus on the design of linear feedback
for the stabilization problem.

For the stability problem we adopt the framework
of polyhedral Lyapunov functions \cite{blanchini_set-theoretic_2015}.
Those are asymmetric norms
(finitely valued, satisfy the triangle inequality,
positively homogeneous, and positive definite)
whose unit balls are polyhedra.
We provide the necessary background on polyhedra and
polyhedral Lyapunov functions in \cref{sec:prelims},
and show how to use them for verifying stability in \cref{sec:stability_conditions}.

The polyhedral Lyapunov framework is less conservative
than standard quadratic Lyapunov functions. 
But this comes at a price:
whereas finding a quadratic Lyapunov function to assess stability
is a convex problem,
efficiently solved using linear matrix inequalities \cite{boyd_linear_1994},
the same problem is hard to solve for polyhedral Lyapunov functions.
In this paper we tackle the problem by developing an approach 
based on a two step iterated procedure,
along the lines of \cite{kousoulidis_optimization_2020}.
Each step is a linear program (LP) and
the iterations build a sequence of polyhedral sets,
each reducing the gap to a contractive set,
i.e. to a decaying Lyapunov function.
The full algorithm is presented in \cref{sec:algorithm}.

A relevant feature of our algorithm is that we fix the complexity of the 
polyhedral Lyapunov function, that is,
the number of vertices or constraints that define the polyhedral set.
This is one of the features that distinguish our approach from 
others available in the literature
\cite{ambrosino_convex_2012,blanchini_set-theoretic_2015,guglielmi_polytope_2017-1,miani_maxis-g:_2005,polanski_absolute_2000-1}.
We provide a comparison in \cref{sec:others}.

Another relevant feature is that we can easily extend the framework
to include additional constraints and new variables,
which is of practical significance for control synthesis.
This allows us to develop most of our discussion for autonomous linear systems,
simplifying the exposition.
We then show how the algorithm is generalized 
to control design and nonlinear systems
in \cref{sec:gen}.
We also show how our algorithm can be used to assess 
contraction/incremental stability of a system.
The effectiveness of the algorithm is illustrated on a DC motor example, 
in \cref{sec:example}.

\section{Mathematical Preliminaries}
\label{sec:prelims}

\subsection{Polyhedra}
We use $\cC$ to denote a convex set
and $\interior \cC$ to denote its interior,
which we define as the largest open subset of $\mathbb{R}^n$ in $\cC$.
Inequalities on vectors and matrices are used in the element-wise sense.
We represent polyhedra in two ways:
\begin{LaTeXdescription}
    \item[V-representation] given $n \times \numrays$ matrix $V$,
    \begin{equation}
        \cCV(V) = \{x: x = V\posvector,\, \posvector \geq 0,\, \innerproduct{1}{\posvector} = 1\}
        \label{eq:V-rep}
    \end{equation}
    \item[H-representation] given $\numrays \times n$ matrix $H$,
    \begin{equation}
        \cCH(H) = \{x: Hx \leq 1 \}
        \label{eq:H-rep}
    \end{equation}
\end{LaTeXdescription}

For this paper we exclusively use V-representation
polyhedra.
All results can be readily extended to 
H-representation polyhedra by duality,
\cite[Proposition 4.35]{blanchini_set-theoretic_2015}.
Geometrically,
the columns of the $V$ correspond to the vertices of $\cCV(V)$,
which is formed by taking their convex hull.

We make the \emph{standing assumptions} that our polyhedra
\begin{align*}
&\bullet \text{ are Bounded} \tag{A1} \label{assump:bounded} \\
&\bullet \text{ contain $0$ in their interior} \tag{A2} \label{assump:absorbing}
\end{align*}
Assumption \eqref{assump:bounded} is always true
for V-representation polyhedra.
Eq. \eqref{assump:absorbing} requires that there exists a $p_+ > 0$
such that $0 = Vp_+$,
which we can test using linear programming (LP).

\subsection{Polyhedral Lyapunov Functions}

Any convex set $\cC$ defines a non-negative function
through its Minkowski functional.
The Minkowski functional of a set $\cC$ is given by 
\begin{equation}
\Psi_\cC(x) = \inf\{r \in \mathbb{R}: r > 0 \text{ and } x \in r\cC \}
\label{eq:minkowski}
\end{equation}

For V-representation polyhedra this means
\begin{equation}
    \psiV(x)
    = \min\{1^T p \,:\, x = Vp,\, p \geq 0 \}
    \label{eq:psiV_primal}
\end{equation}
Using the strong duality of LP problems \cite[Section 5.2.1]{boyd_convex_2004},
this is equivalent to 
\begin{equation}
    \psiV(x)
    = \max\{h^T x \,:\, h^T V \leq 1^T \}
    \label{eq:psiV_dual}
\end{equation}
Using LP, \eqref{eq:psiV_primal} and \eqref{eq:psiV_dual} can be evaluated
numerically at any given point.
For V-representation polyhedra
satisfying \eqref{assump:absorbing},
$\psiV$ is always a real valued function, strictly positive for all $x \neq 0$,
and radially unbounded.
As such they are valid \emph{Lyapunov function candidates}.

However, for us to fully establish the link with Lyapunov theory
we also require the notion of a derivative of the
candidate function.
Unfortunately, the Minkowski functionals of polyhedra are
not everywhere differentiable and we need to resort
to the more general setting of subdifferentials:
\begin{equation}
    \partial\Psi(x) = \{h \,:\, \forall y\in \mathbb{R}^n, \, \Psi(y)-\Psi(x) \geq h^T(y-x) \}
    \label{eq:subdiff}
\end{equation}
For differentiable functions the subdifferential at $x$ corresponds to the gradient.

For Minkowski functionals of V-representation polyhedra,
we can obtain an explicit representations of $\partial\Psi_V(x)$.
\begin{proposition}[V-representation Subdifferential]
    For any given matrix $V$, consider the Minkowski functional $\Psi_V$.
    Then, 
    \begin{equation}
        \partial\psiV(x) = \{h \,:\, h^T x = \psiV(x),\,h^TV \leq 1^T \}
        \label{eq:subdiff_v}
    \end{equation}
    \label{prop:subdiff_v}
\end{proposition}
\begin{proof}
    $\eqref{eq:subdiff_v} \subseteq \partial\psiV(x)$:
    We need to show that
    $\psiV(y) \geq h^T y$.
    From \eqref{eq:psiV_dual}, this must be true.

   $\partial\psiV(x) \subseteq \eqref{eq:subdiff_v}$:
   We prove this in two steps.
   We first show that any $h \in \partial\psiV(x)$ must
   satisfy $h^T V \leq 1^T$ by contradiction.
   If $h$ doesn't satisfy $h^T V \leq 1^T$,
   there must exist a column of $V$, $[V]_i$,
   such that $h^T [V]_i = c > 1$ and $\psiV([V]_i) = 1$.
   If we set $y = k [V]_i$, for $h \in \partial\psiV(x)$
   we must then have $k - \psiV(x) \geq kc - h^T x$
   for all $k > 0$.
   However, if we set
   $k > \max(0,(h^Tx - \psiV(x))/(c-1))$
   the inequality is invalidated, leading to a contradiction.
   As such, if $h \in \partial\psiV(x)$ then $h^T V \leq 1^T$.
   If we set $y = 0$ we get
   $h^Tx \geq \psiV(x)$.
   But from before $h^T V \leq 1^T$, and so $h$ must be
   included in the domain of \eqref{eq:psiV_dual} and
   therefore $h^T x = \psiV(x)$.
   Since $h$ satisfies both
   $h^T V \leq 1^T$ and $h^T x = \psiV(x)$,
   $h \in \eqref{eq:subdiff_v}$.
\end{proof}

We use subdifferentials to quantify how our candidate Lyapunov function $\psiV$
changes along system trajectories.
\begin{proposition}
    For any trajectory $x(\cdot)$ of the autonomous system $\dot{x} = Ax$,
    \begin{equation}
        \psiV(x(t)) - \psiV(x(0)) = \int_0^t
        \sup_{h \in \partial \psiV(x(\tau))} h^T Ax(\tau) d\tau
    \label{eq:subdiff_derivative}
    \end{equation}
\end{proposition}
\begin{proof}
    See \cite{queiroz_lyapunov_2004} and
    \cite[Sections 2.2]{blanchini_set-theoretic_2015}.
\end{proof}

The integrand is equivalent to
$\max h^T Ax$ subject to $h^T x = \psiV(x),\,h^TV \leq 1^T$,
which is an LP that we can directly numerically compute.

\section{Contraction, Linear Programming, \\ and Lyapunov Stability}
\label{sec:stability_conditions}

Convex sets can be used to show stability through
the notion of contraction
(also frequently referred to as strict positive invariance).
For continuous time systems,
this requires us to first define
tangent cones.
For a general convex set $\cC$,
the interior of the tangent cone
at a point $x \in \cC$ is given by
\begin{equation}
    \tconeint{\cC}(x) = \{ w \,:\,
    \lambda > 0,\,
    x + \lambda w \in \interior \cC \}
    \label{eq:tconeint}
\end{equation}

We then say that $\cC$ is a \emph{contracting convex set
under the dynamics} $\dot{x} = Ax$ if
\begin{equation}
    \forall x \in \cC ,\, Ax \in \tconeint{\cC}(x)
    \label{eq:geo_contraction}
\end{equation}

This can be interpreted in words as saying that 
the flow induced by the dynamics maps all points in $\cC$
to points in $\interior \cC$ for all $t > 0$.
For polyhedra,
we can express this geometric condition in a more convenient form.

\begin{proposition}[Contraction in V-representation]
    \label{prop:contracting-v}
    For a polyhedron $\cCV(V)$,
    \eqref{eq:geo_contraction}
    is equivalent to
    $\exists 
    M \in \mathbb{R}^{\numrays \times \numrays}$
    and scalar $\rate > 0$
    such that
    \begin{subequations}
        \label{eq:contracting}
        \begin{align}
            -\rate 1^T &= 1^T M, \label{eq:contracting_rate} \\
            AV &= VM,
            \label{eq:contracting_equality_main} \\
            \offdiag(M) &\geq 0 \text{ [$\equiv M$ is Metzler]}\label{eq:contracting_metzler}
        \end{align}
    \end{subequations}

    Where we use $\offdiag(M) \geq 0$ to denote that all elements
    of $M$ outside of the main diagonal are non-negative.
    We call $\rate$ the rate of contraction.
\end{proposition}
\begin{proof}
    For $\cCV(V)$, we use \eqref{eq:V-rep} to express \eqref{eq:geo_contraction} as
    \begin{alignat*}{1}
        \forall &(p \geq 0, 1^Tp = 1): \\
        \exists &(\lambda > 0,\hat{p}\geq 0,1^T\hat{p} < 1):
        AVp = \lambda^{-1} V(\hat{p} - p)
    \end{alignat*}
    This is equivalent to
    \begin{equation*}
        \exists (\Lambda > 0, \hat{P} \geq 0, 1^T \hat{P} < 1^T):
        AV = V\Lambda^{-1}(\hat{P}-I)
    \end{equation*}
    where $\Lambda$ is a diagonal matrix.
    Setting $\hat{M} = \Lambda^{-1}(\hat{P}-I)$,
    we rewrite this as:
    \begin{equation*}
        \exists (\rate > 0, \hat{M}):
        1^T \hat{M} \leq -\rate 1^T, AV = V\hat{M}, \offdiag(\hat{M}) \geq 0
    \end{equation*}
    From \eqref{assump:absorbing}, $\exists p_+ > 0$ such that $0 = Vp_+$.
    We can then add this $p_+$ to each column of $\hat{M}$ until it satisfies \eqref{eq:contracting_rate}
    without affecting \eqref{eq:contracting_equality_main} or \eqref{eq:contracting_metzler}.
\end{proof}

Testing contraction is thus an LP optimization problem.
\begin{alignat}{2}
   &\max_{M,\rate}\, &\> &\rate
   \label{eq:contracting_optimization} \\
   &\text{s.t.:}&
        &\eqref{eq:contracting} \nonumber
\end{alignat}
The conditions of \cref{prop:contracting-v} are verified
if the maximization in \eqref{eq:contracting_optimization} is greater than $0$.

It is fairly straightforward to show that set contraction implies stability
through the homogeneity of the linear dynamics.
However, using \eqref{eq:subdiff_v},
we can build a direct connection to Lyapunov theory.
\begin{proposition}[Lyapunov stability]
    $\cCV(V)$ is a contracting set for $A$ with rate $\rate$
    if and only if
    \begin{equation}
        \forall x \in \mathbb{R}^n \quad 
        \sup_{h \in \partial \psiV(x)} h^T Ax
        \leq -\rate \psiV(x)
        \label{eq:lyap_contracting}
    \end{equation}
    \label{thm:contraction_equiv}
\end{proposition}
\begin{proof}
    \eqref{eq:contracting} $\implies$ \eqref{eq:lyap_contracting}:
    Starting with \eqref{eq:subdiff_v}
    and \eqref{eq:psiV_primal}
    we have that there exists some $p$ such that
    $h^Tx = 1^Tp = \psiV(x),\, x = Vp,\, p \geq 0$
    for the $h \in \partial \psiV(x)$ that maximizes $h^T Ax$.
    Also note that we can rewrite a matrix $M$
    satisfying \eqref{eq:contracting_metzler}
    as $M = P + \delta I$,
    where $P$ is an positive matrix.
    Therefore
    \begin{align*}
        h^TAx
        &= h^TAVp \\
        &= h^TVMp \quad \eqref{eq:contracting_equality_main} \\
        &= h^TV(P + \delta I)p = h^TVPp + \delta h^TVp \\
        &= h^TVPp + \delta 1^Tp \quad [x = Vp,\, 1^Tp = h^Tx] \\
        &\leq 1^TPp + \delta 1^Tp = 1^TMp \quad
        \text{[$h \in \partial \psiV(x),\, Pp \geq 0$]} \\
        &\leq -\rate 1^T p \quad \eqref{eq:contracting_rate} \\
        &\leq -\rate \psiV(x)
    \end{align*}

    \eqref{eq:lyap_contracting} $\implies$ \eqref{eq:contracting}:
    We show this by contradiction.
    The optimization problem in 
    \eqref{eq:contracting} can be broken down to multiple
    multiple sub-problems, one for each column of $V$, $[V]_i$:
    \begin{equation}
        \exists p \in \mathbb{R}^{m} \,:\,
        \rate = -1^Tp,\,
        A[V]_i = Vp,\,
        p_{j \neq i} \geq 0
    \label{eq:contracting_individual}
    \end{equation}
    where $\{j \neq i\}$ is $[1,m] \setminus \{i\}$.
    If \eqref{eq:contracting} is infeasible for a given $\rate$
    then \eqref{eq:contracting_individual} must be infeasible
    for at least one $[V]_i$.
    This is because otherwise a feasible solution
    for \eqref{eq:contracting} could be generated by
    aggregating $\{p\}$
    into a matrix $P$.
    As such we begin by assuming that there exists some $[V]_i$
    for which \eqref{eq:contracting_individual} is infeasible.
    Then, by Farkas' Lemma, there exists an $h$ such that
    \begin{equation}
        h^TA[V]_i > -\rate ,\,
        h^T[V]_i = 1 ,\,
        h^TV \leq 1^T
    \end{equation}
    This means that $h \in \partial \psiV([V]_i)$.
    But, since, 
    $\psiV([V]_i) = 1$,
    $h^TA[V]_i > -\rate \psiV([V]_i)$,
    contradicting \eqref{eq:lyap_contracting}.
\end{proof}

From Lyapunov theory, 
\eqref{eq:lyap_contracting} with $\rate > 0$ combined 
to \eqref{eq:subdiff_derivative} imply global exponential stability of the system
with convergence rate $\rate$.
The connection between contraction and 
Lyapunov stability is a well-known result in system theory
\cite[Thm 4.24, 4.33]{blanchini_set-theoretic_2015}.
However, we include this proof to show how we can use
\eqref{eq:subdiff_v}
to directly proof results about polyhedral Lyapunov functions.
We hope to use this approach to extend results that have no
direct interprettation as set contraction,
such as dissipativity,
to polyhedral Lyapunov functions.
This can be particularly useful for control design and for the analysis
of open systems.

\section{Algorithm for Finding Polyhedra}
\label{sec:algorithm}

\subsection{Overview}
Testing if a set is contracting under the linear dynamics $\dot{x} = Ax$
is a much easier problem than finding a contractive set.
In fact, treating 
$V$ as a variable makes 
\eqref{eq:contracting_equality_main}  nonconvex and nonlinear
(since $M$ and $\rate$ are also variables).
However, we can instead try to find such a $V$ iteratively,
by decomposing \eqref{eq:contracting_optimization} into a pair of LP problems.
One considers $V$ fixed and estimates \emph{the contraction gap of $V$}.
The other,  \emph{introduces variations on $V$ to reduce the contraction gap}.

To initialize our algorithm we need to build an initial
candidate polyhedron which satisfies our standing assumptions.
A possible construction is suggested in Remark \ref{rem:init}.
The number of vertices $\numrays$ of our polyhedron will remain the same
throughout the whole procedure.
This means that we can control the \emph{complexity} of our Lyapunov function
and treat $\numrays$ as an input parameter.

\textbf{(i)~Estimation of the contraction gap:}
for any given candidate polyhedron $V$, the contraction gap is estimated
via the LP optimization  \eqref{eq:contracting_optimization}
The following Proposition justifies the use of \eqref{eq:contracting_optimization}.
\begin{proposition}
    Given a $\cCV(V)$ that satisfies \eqref{assump:absorbing},
    there exists some finite $\rate$ (not necessarily $ > 0$) that satisfies
    \eqref{eq:contracting}.
    Additionally, if $\cCV(V)$ satisfies \eqref{eq:contracting}
    for some $\rate^*$,
    it will also satisfy \eqref{eq:contracting} for any $\rate < \rate^*$.
    \label{proposition:rate_as_distance}
\end{proposition}
\begin{proof}
    Note that \eqref{assump:absorbing} implies that $V$
    is full row rank since otherwise the strict interior of
    $\cCV(V)$ would be empty.
    We rewrite $p = p_{b} + \delta_i p_+$,
    where $p_{b}$ is a solution to 
    $A[V]_i = Vp_{b}$
    (implied to exist because $V$ is full row rank)
    and $p_+ > 0$ is a positive solution to $0 = Vp_+$
    (implied to exist because of \eqref{assump:absorbing}).
    There then exists a $\delta_i^*$ such that
    $A[V]_i = Vp$ and $p_{i\neq j} \geq 0$ for all
    $\delta_i \geq \delta_i^*$.
    Since $1^Tp_+ > 0$, increasing $\delta_i$
    decreases $\rate_i = -1^Tp_{b}-\delta_i 1^Tp_+$.
    As such, there exists $\{ \delta_i \}$ such that
    for all $i$ 
    $\rate_i^* = -\min_i\{1^Tp_{b} + \delta_i^*p_+\}$.
    Aggregating $\{p\}$ into a matrix $P$
    produces a solution to \eqref{eq:contracting} with
    $\rate = \rate_i^*$.
    By increasing $\{\delta_i\}$ further,
    \eqref{eq:contracting} can be satisfied for any
    $\rate < \rate^*$.
\end{proof}

If $\rate > 0$ we conclude that our  polyhedron
is contracting and we are done. 

\textbf{(ii)~Reduction of the contraction gap:}
if $\rate \leq 0$, the following LP program searches for 
variations on $V$ that are compatible with a reduction of
the contraction gap.
For a small $\varepsilon > 0$, solve:
    \begin{subequations}
        \label{eq:LP2}
        \begin{align}
           \max_{\delta \eta , \delta M, \delta V} \, &\> \delta \rate \\
            -\delta \rate 1^T &= 1^T \delta M, \label{eq:LP2_contracting_rate} \\
            A \delta V &= \delta V M + V \delta M,
            \label{eq:LP2_equality_main} \\
            0 &\leq \offdiag(M+\delta M) \\
   \norm{\delta V}{1} & \leq \varepsilon ,\,
   \norm{\delta M}{1} \leq \varepsilon 
        \end{align}
    \end{subequations}
The idea is to look for variations $\delta V$ in a small neighborhood of $V$
that would improve the gap $\eta < \eta + \delta \eta$.
The search space is limited to the tangent space of the constraints manifold,
at the point $(\eta, V, M)$.
Therefore, for small enough variations, any solution to the LP above improves the contraction rate 
while preserving the integrity of the constraints.
We thus build the new polyhedron as $V_{\textrm{new}} := V + \delta V$.
A small $\varepsilon$ ensures that the $V_{\textrm{new}}$ remains in a small neighborhood of  $ V $
to be compliant with the linear approximation above.
It also guarantees that $\cCV(V_{\textrm{new}})$ satisfies our standing assumption (by continuity argument). 

The procedure iterates between these two steps until the estimation of the contraction gap
returns a positive $\eta$.

\begin{remark}
\label{rem:init}
We compute
an initial candidate polyhedron
by randomly sampling $\numrays-1$ points
in $\mathbb{R}^n$
and rescaling them to have length 1.
We then add an additional vertex equal to minus their sum
and again normalised to have length 1.
This way we obtain 
a polyhedron $\cCV(V)$
with $\numrays$ irreducible vertices 
that, as long as $\numrays \geq n + 1$,
must also satisfy \eqref{assump:absorbing}.
\end{remark}

\subsection{Fast reduction of the contraction gap}
\label{sec:fast_gap}

We can speed up \eqref{eq:LP2} by not having $\delta \rate$ and $\delta M$
as variables in the second optimization problem.
We do this 
by considering directions in the tangent space
that are compatible with the optimality conditions of \eqref{eq:contracting_optimization}
around a fixed previous solution.
With these extra restrictions in place, we can obtain a directional derivative
of $\rate$ with respect to $V$.

To estimate this derivative we need to rewrite \eqref{eq:contracting_optimization}
into standard LP form 
\begin{subequations}
\label{eq:contracting_optimization_sf}
    \begin{alignat}{2}
       &\max_{\hat{x}}\, &\> 
        \hat{c}^T\hat{x}&
       \label{eq:contracting_optimization_sf_obj}
       \\
       &\text{s.t.:}&
       \hat{A}\hat{x} &= \hat{b}
        \label{eq:contracting_optimization_sf_eq}
        \\
        & &   \hat{x} &\geq 0
        \label{eq:contracting_optimization_sf_ineq}
    \end{alignat}
\end{subequations}
where $\hat{A}$ is a $k \times l$ matrix.
A detailed transformation procedure is provided in \cref{subsec:grad_detail}.

Using the standard form, we can then make use of the notion of Basic Feasible Solutions (BFS):
from \cite[Section 13.2]{nocedal_numerical_2006},
a BFS of \eqref{eq:contracting_optimization_sf}
is a set $\cB \subset \{1,.., l\}$
and a vector $\hat{x}$
such that:
\begin{itemize}
    \item $\hat{A}\hat{x} = \hat{b}$ and $\hat{x} \geq 0$
    \item $\cB$ contains $k$ elements
    \item $i \not\in \cB \implies \hat{x}_i = 0$.
    \item The $k \times k$ matrix $\hat{B}$ defined as
    \begin{equation*}
        \hat{B} = [\hat{A}_i]_{i \in \cB}
    \end{equation*}
    is non-singular, where $\hat{A}_i$ is the $i$th column of $\hat{A}$.
\end{itemize}
When solving \eqref{eq:contracting_optimization_sf}
it is sufficient to only consider BFSs
\cite[Theorem 13.2]{nocedal_numerical_2006}.
We can obtain a BFS 
that is also
an optimal solution of \eqref{eq:contracting_optimization_sf}
via the Simplex algorithm \cite{nash_origins_1990}.
Many commercial interior point solvers
also offer the option of deriving a BFS
as a cheap post-processing stage.

With the optimal BFS we can now find the derivative of $\rate$
with respect to changes in $V$. 
In general, we need to implicitly differentiate
equation \eqref{eq:contracting_optimization_sf_eq}
to obtain a derivative of $\hat{x}$
(which represents our variables $\rate$ and $M$ in standard form)
with respect to $\hat{A}$ and $\hat{b}$
(which are functions of $V$).
However, an explicit solution can only be obtained if $\hat{A}$ is invertible
\cite{dontchev_implicit_2009}.
Instead, 
we restrict ourselves to the BFS and consider only
variables that aren't constrained to be $0$ by our choice of $\cB$.
Indeed, given $\hat{y} = \{\hat{x}_i : i \in \cB\}$,
we can always obtain the derivative of $\hat{y}$
with respect to $\hat{B}$ and $\hat{b}$. 
Because $\hat{B}$ is invertible
(from the properties of the BFS),
this derivative is always well defined.

Using the transformation in \cref{subsec:grad_detail},
we can map the LP problem back to its original coordinates
and build an explicit representation of 
$D_V \rate$ and $D_V M$.
We can then replace \eqref{eq:LP2} with
the following LP optimization problem:
\begin{alignat}{2}
   &\max_{\delta V}\, &\> &(D_V \rate) (\delta V)
   \label{eq:second_step} \\
   &\text{subject to: }& &\norm{\delta V}{1} \leq \varepsilon 
    \nonumber \\
   &\text{and}:&\quad
        \offdiag(M + &(D_V M) (\delta V)) \geq 0
    \nonumber
\end{alignat}
By constraining $\delta V$ to be small enough
we can also ensure that the new polyhedron
$\cCV(V+\delta V)$ satisfies
\eqref{assump:bounded} and \eqref{assump:absorbing}.
This can also be enforced directly,
as in \cite{kousoulidis_optimization_2020},
but keeping $\delta V$ small
seems to work sufficiently well in practice.

We summarize the full algorithm below:
\begin{algorithm}[H]
  {\bf Data:} The matrix $A$, the maximum number of iterations, \\
  and the number of vertices in the matrix $\numrays$ \\
  {\bf Result:} $V$ satisfying \eqref{eq:contracting} if found, else $False$ \smallskip \\
  {\bf Procedure:} \\
  $V^{(0)} = $ Initialize$(m)$ \\
  $k = 0$ \\
  {\bf while }{$k \leq \texttt{max\_iter}$}: \\
  \hspace*{4mm} $\rate, M, \cB = $ 
  Optimal BFS of \eqref{eq:contracting_optimization}
  with $V = V^{(k)}$ \\
  \hspace*{4mm} {\bf if }{$\rate > 0$}: \\
  \hspace*{8mm} return $V^{(k)}$ \\
  \hspace*{4mm} $D_V \rate, D_V M \!=\! $
  Derivatives of \eqref{eq:contracting_optimization}
  using $\cB$ and $V \!=\! V^{(k)}$ \\
  \hspace*{4mm} $\delta V = $ \eqref{eq:second_step} using
    $M,\, D_V \rate,$ and $D_V M$
   \\
  \hspace*{4mm} $V^{(k+1)} = V^{(k)} + \delta V$ \\
  \hspace*{4mm} $k := k + 1$ \\
  return $False$ \\
  \caption{Stability Verification Algorithm \label{alg:verification}}
\end{algorithm}

\subsection{The Derivative Computation in Detail}
\label{subsec:grad_detail}

To go from \eqref{eq:contracting_optimization}
to \eqref{eq:contracting_optimization_sf}
we use the vectorization operator $\vect(\cdot)$,
which converts
a $n \times m$ matrix into
a column vector with $nm$ entries
by stacking up its columns
and the Kronecker product ($\otimes$):
\begin{alignat}{2}
    \hat{A} =& 
             \begin{bmatrix}
                I_m \otimes V & -\vect(V) & 0_{nm} & 0_{nm} \\
                I_m \otimes 1_m^T & -1_m & 1_m & -1_m
            \end{bmatrix}
    &
    \quad\hat{b} =&
       \begin{bmatrix}
           \vect(AV) \\
           0_m
       \end{bmatrix}
       \nonumber
    \\
    \hat{c} =&
            \begin{bmatrix}
                0_{mm} \\
                0 \\
                1 \\
                -1 
            \end{bmatrix}
    &
    \quad\hat{x} =&
       \begin{bmatrix}
           \vect(P) \\
           -\delta \\
           \rate_+ \\
           \rate_-
       \end{bmatrix}
    \label{eq:standard_form_transformations}
\end{alignat}
where $\rate_+ - \rate_- = \rate$ and $P + \delta I = M$.
We can also rewrite
\eqref{eq:contracting_optimization_sf_eq}
as
\begin{equation}
    f = \hat{A}\hat{x} - \hat{b} = \hat{B}\hat{y} - \hat{b} = \hat{G}\vect(V) = 0
    \label{eq:grad_constraints}
\end{equation}
where
\begin{equation}
    \hat{G} = \begin{bmatrix}
        M^T \otimes I_n - I_m \otimes A \\
        0_{m \times (nm)}
    \end{bmatrix}
\end{equation}
From \eqref{eq:grad_constraints},
we use the implicit function theorem to compute
the Jacobian of $\hat{y}$ with respect to $\vect(V)$
\begin{equation}
   \nabla_{\vect(V)} \hat{y}
   = -(\nabla_{\hat{y}} f)^{-1} \nabla_{\vect(V)} f
   = -\hat{B}^{-1}\hat{G}
   \label{eq:jacobian}
\end{equation}
such that $(D_V \hat{y}) (\delta V) = (\nabla_{\vect(V)} \hat{y}) \vect(\delta V)$.

For our derivatives of interest
\begin{alignat*}{1}
    (D_V \rate) (\delta V)  = (D_V \rate_+) (\delta V) &- (D_V \rate_-) (\delta V) \\
    \text{and, }\forall i \neq j,\, (D_V M_{ij}) (\delta V) &= (D_V P_{ij}) (\delta V)
\end{alignat*}
Note that all derivatives on the right hand side are of elements of $\hat{x}$.
To compute each of them
we first obtain the corresponding index $i$ in $\hat{x}$
from \eqref{eq:standard_form_transformations}.
Then,
if $i \not\in \cB$, the derivative is zero.
Else, we obtain the corresponding
index $j$ in $\hat{y}$ so that $\hat{x}_i = \hat{y}_j$.
To compute the derivative,
we then multiply
the $j^{th}$ row of \eqref{eq:jacobian}
with $\vect(\delta V)$.

Our derivative calculations can be interpreted as only considering the effect of changes
in $V$ to our variables under a fixed basis.
As long as the solution with the current basis is feasible for $V + \delta V$,
this will still be a BFS 
(although it might not be the optimal one).

\section{Generalizations}
\label{sec:gen}

\subsection{Feedback Design}
\label{sec:feedback}

\cref{alg:verification} can be readily extended for output-feedback design for systems
of the form
$$
\dot{x} = A x + B u \qquad\quad y = C x.
$$
This includes state feedback as a special case
when $C = I$.
We replace \eqref{eq:contracting_optimization} with
\begin{alignat}{2}
   &\max_{M, K, \rate}\, &\> &\rate
   \label{eq:contracting_optimization_design} \\
   &\text{s.t.:}&
        -\rate 1^T &= 1^T M, \nonumber \\
            && (A + BKC)V &= VM,
            \nonumber \\
            && \offdiag(M) &\geq 0 \text{ [$\equiv M$ is Metzler]}
            \nonumber
\end{alignat}

The two step procedure outlined in Section \ref{sec:algorithm} 
reduces the nonlinear (for non-fixed $V$) program \eqref{eq:contracting_optimization}
into an iteration involving linear programs.
A similar reduction can be 
applied to \eqref{eq:contracting_optimization_design}.
Just like \eqref{eq:contracting_optimization},
\eqref{eq:contracting_optimization_design} is a linear program for 
fixed $V$, which means that the estimation of the contraction 
gap is a straightforward operation.
Likewise, for \eqref{eq:contracting_optimization_design},
the reduction of contraction gap step is performed on a 
larger manifold, which now also contains $K$.
However, the search procedure remains essentially the same.
This means that the approach in \eqref{sec:fast_gap}
can also be adapted to control synthesis (it just requires different transformations to 
get the constraints into standard form). 

We can also consider scenarios beyond simple stabilization.
New formulations may be considered 
as long as they are encoded into linear constraints.
Natural directions to explore are testing/enforcing input/output gains,
the design of dynamic output feedback,
and the design of nonlinear state-feedback controllers taking advantage of the Lyapunov framework.
This will be the goal of future research.

\subsection{From Linear Stability to LDIs}
\label{sec:ldi}

The discussion in this paper has been developed for simple linear
systems but,
since \eqref{eq:contracting_optimization} and \eqref{eq:contracting_optimization_design} are convex,
the same approach can be used for the analysis LTV and LPV systems, and even for nonlinear systems.

The idea is that systems of the form
\begin{equation}
    \dot{x} = A(w(t))x + B(w(t))u
    \label{eq:polytopic_ldi}
\end{equation}
satisfy, for all $t \geq 0$,
\begin{equation*}
 A(w(t)) \in \convexhull\{A_i\},\, B(w(t)) \in \convexhull\{B_i\}
\end{equation*}
where
$$\convexhull\{A_i\} = \{A: A = \sum_i p_i A_i,\, p_i \geq 0,\, \sum_i p_i = 1 \};$$
and likewise for $\convexhull\{B_i\}$.

Stability is thus verified as in \eqref{eq:contracting_optimization}
where constraints \eqref{eq:contracting_rate} and \eqref{eq:contracting_equality_main} are 
replaced by several constraints of the form
$$
-\rate 1^T = 1^T M_i \text{ and } A_iV = VM_i \ .
$$
Indeed, the extension of the constraint set does not change the structure
of the solution. Combined with \eqref{eq:contracting_optimization_design},
this allows us to tackle the design of a stabilizing linear state-feedback law $u = Kx$.

A similar approach can be explored to verify contraction
\cite{lohmiller_contraction_1998-1, Russo2010, Pavlov2005, forni_differential_2014} 
of nonlinear systems of the
form  $\dot{x} = f(x)$.
The problem reduces to finding a contracting polyhedral set $V$ for the linearization
$$
\dot{\delta x} = \partial f(x) \delta x \qquad \quad \forall x \ .
$$
The solution is  obtained as above, 
by relaxing the problem to a suitable convex hull  $\partial f(x) \in \convexhull\{A_i\}$
for all $x$. 

\section{Comparison with Related Approaches}
\label{sec:others}

In systems theory polyhedral Lyapunov functions are usually constructed using
algorithms that progressively add new vertices or constraints until
a contractive set is reached.
Other algorithms exploit optimization based approaches to directly construct
contractive sets.

\cite{blanchini_set-theoretic_2015}, \cite{guglielmi_polytope_2017-1} and \cite{miani_maxis-g:_2005}
provide examples of the iterative approach.
The new vertices are generated by iterating the update equation
of discrete-time linear systems.
The main advantage of these approaches is that in many cases the
sequence of polyhedra produced can be guaranteed to converge to
a polyhedral Lyapunov function, if one exists.
The main disadvantage is that not all linear constraints can be
incorporated easily.
Additionally,
it is difficult to bound the complexity
of the resulting polyhedron.

Examples of optimization-based algorithms can be found in
\cite{ambrosino_convex_2012} and \cite{polanski_absolute_2000-1}.
\cite{polanski_absolute_2000-1}
considers positive rescalings
of the vertices of the initial candidate polyhedron.
This allows for fast computation
but may lead to high complexity polyhedra
even when a simple contractive polyhedron exists.
The algorithm 
in \cite{ambrosino_convex_2012}
makes clever use of both representations of a polyhedron.
However, this is also a limiting factor
because of the possible difference in the number
of vertices and constraints of the two representations.
In higher dimensions,
a few vertices/constraints in one representation
may lead to a very large number of constraints/vertices
in the other representation \cite{mcmullen_maximum_1970},
which hinders the use of the algorithm in high dimensional settings.

The approach proposed in this paper is also optimization-based.
It takes into account explicitly the problem of complexity of the polyhedron
by taking the number of vertices of the polyhedron
as a parameter.
This comes at the price of weaker convergence guarantees
(if the predefined number of vertices is too small the problem is unfeasible).
Our algorithm generalizes the approach of \cite{polanski_absolute_2000-1}
by allowing vertices to move freely,
and, in constrast to \cite{ambrosino_convex_2012},
only uses a single representation.
Our approach also tackles the problem linear control synthesis
which is not considered by any of the approaches above.

\section{Example: DC Motor}
\label{sec:example}

We test our optimization-based approach
on the DC motor speed and position examples from
\cite{ctms_control_tutorials}.
For the DC motor speed example,
we carry out a stability analysis under
various possible parameter variations.
We then solve the feedback design problem for the DC motor position example
using both static and dynamic state feedback.
We use the parameters given in the website for the speed example
as our nominal values throughout
(inertia $J_0=.01$, viscous friction constant $b_0=.1$, emf constant $K_0=.01$, resistance $R_0=1$, and inductance $L_0=.5$).

\begin{figure}[tbp]
   \begin{minipage}{.49\linewidth}
        \includegraphics[width=\linewidth]{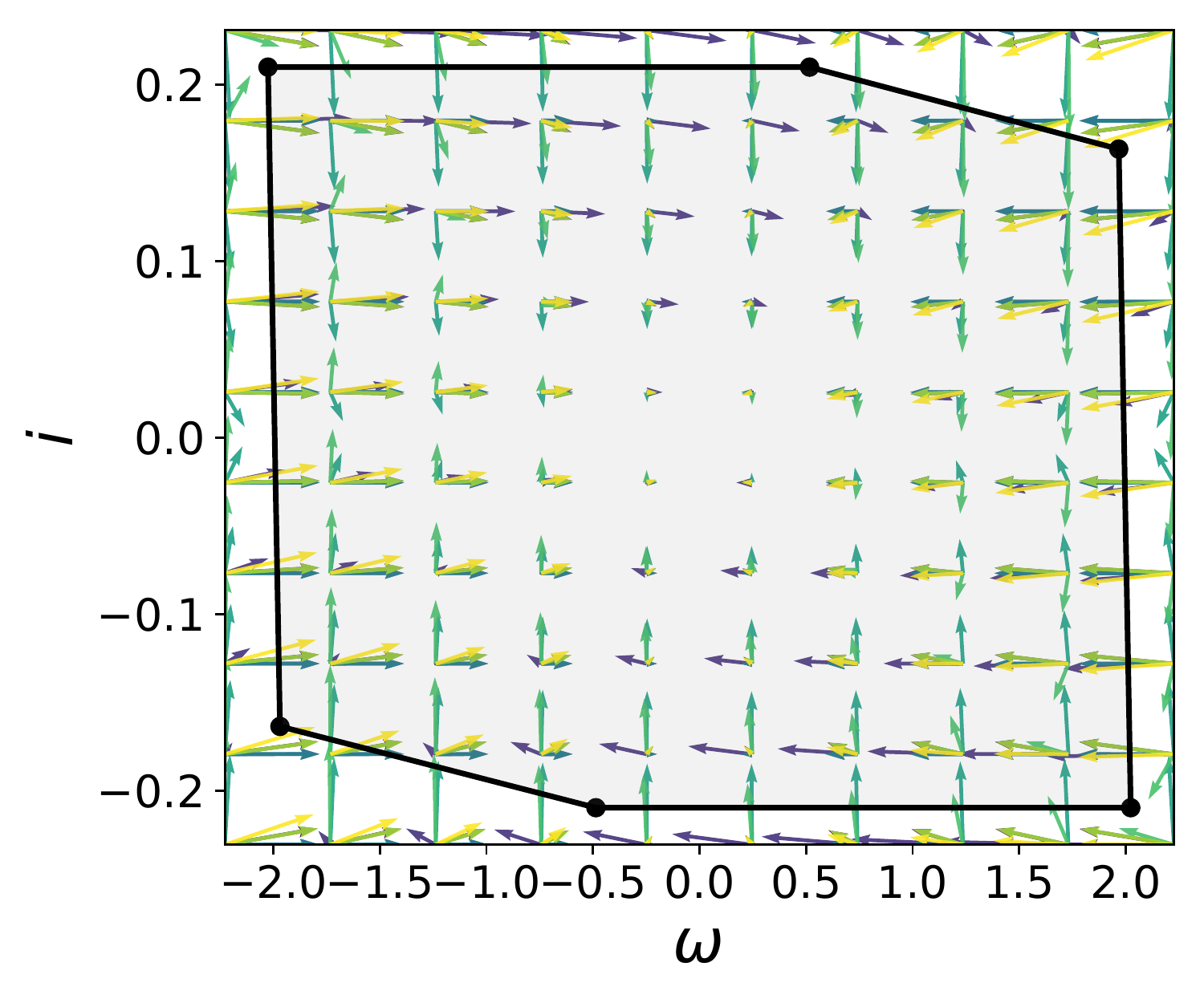}
   \end{minipage}
   \begin{minipage}{.49\linewidth}
        \includegraphics[width=\linewidth]{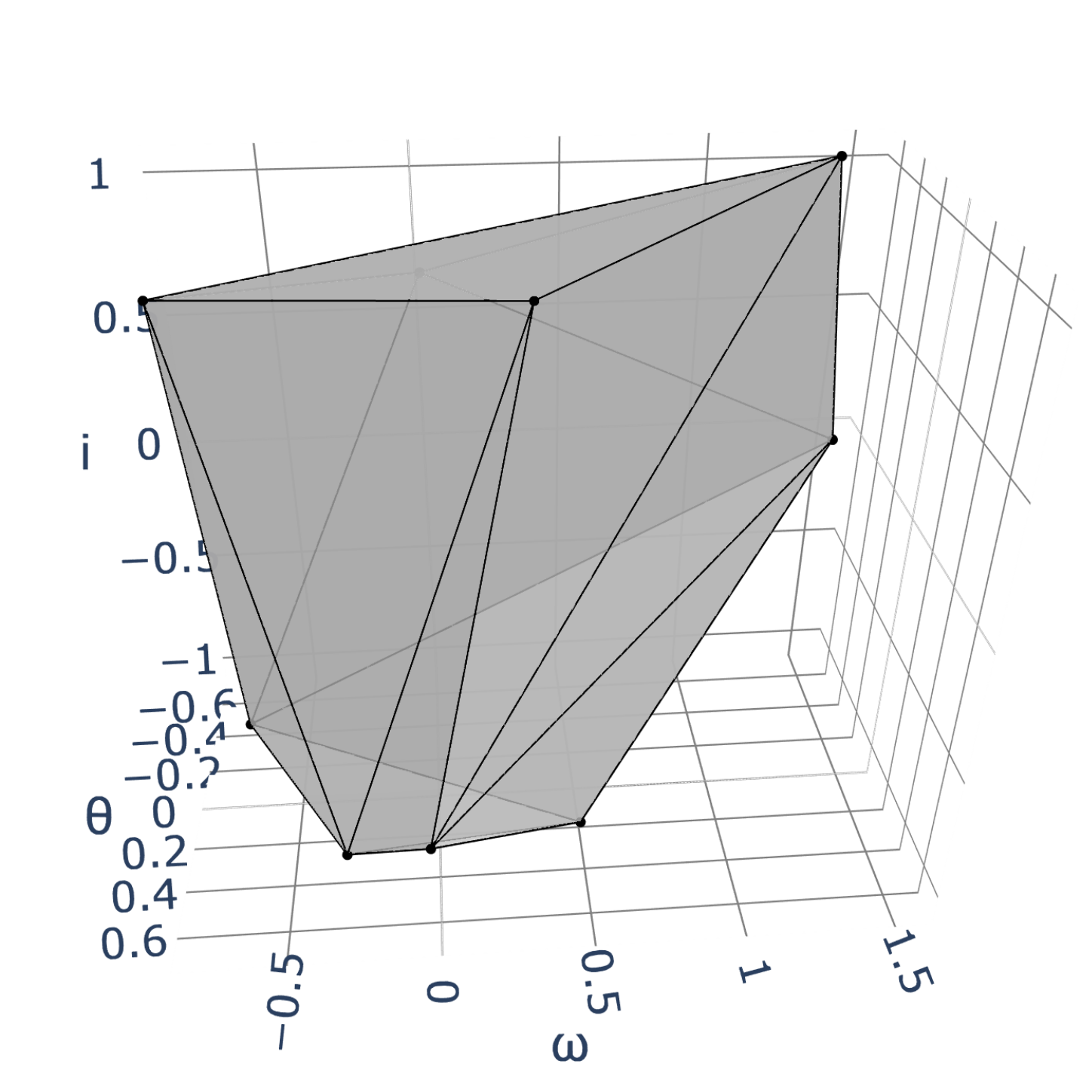}
   \end{minipage}
   \caption{
   \textbf{Left:} Invariant polyhedron for \eqref{eq:motor_speed} and \eqref{eq:motor_speed_params}, \\
    for $\gamma_s = 10$
    ($\rate \geq 0.07$). 
   \textbf{Right:} Invariant polyhedron for \eqref{eq:motor_position} and \eqref{eq:motor_position_params}, 
   with $\gamma_p = 4$ and $u$ as in \eqref{eq:motor_position_feedback}
   ($\rate \geq 0.004$, $0.03$ for nominal)
   }
   \label{fig:motor}
 \end{figure}

The DC motor speed example is governed by the following dynamics:
\begin{equation}
    \frac{d}{dt}
    \begin{bmatrix}
        \omega \\
        i
    \end{bmatrix} 
    = \begin{bmatrix}
        -b/J & K/J \\
        -K/L & -R/L
    \end{bmatrix}
    \begin{bmatrix}
        \omega \\
        i
    \end{bmatrix} 
    \label{eq:motor_speed}
\end{equation}
In this case,
we consider 2 orders of magnitude of parameter perturbations of the form
\begin{equation}
J \in [J_0/\gamma_s,\gamma_s J_0],\, b \in [b_0/\gamma_s,\gamma_s b_0],\, K \in [K_0/\gamma_s,\gamma_s K_0] 
\label{eq:motor_speed_params}
\end{equation}
for $\gamma_s = 10$.
Using \cref{alg:verification},
we are able to find a polyhedron
satisfying \eqref{eq:lyap_contracting}
and \eqref{eq:motor_speed_params} (with $\rate = .07$).
It is shown in \cref{fig:motor}
and has only 6 vertices.
For comparison,
we are only able to find a common quadratic function
using LMIs up to around $\gamma_s = 8.7$.

We next consider the DC motor position example.
The dynamics are now:
\begin{equation}
    \frac{d}{dt}
    \begin{bmatrix}
        \theta \\
        \omega \\
        i
    \end{bmatrix} 
    = \begin{bmatrix}
        0 & 1 & 0 \\
        0 & -b/J & K/J \\
        0 & -K/L & -R/L
    \end{bmatrix}
    \begin{bmatrix}
        \theta \\
        \omega \\
        i
    \end{bmatrix} 
    + \begin{bmatrix}
        0 \\
        0 \\
        1/L
    \end{bmatrix}
    u
    \label{eq:motor_position}
\end{equation}
And we assume that we can only measure $\theta$ and $i$:
\begin{equation}
    y = \begin{bmatrix}
        1 & 0 & 0 \\
        0 & 0 & 1
    \end{bmatrix}
    \begin{bmatrix}
        \theta \\
        \omega \\
        i
    \end{bmatrix} 
\end{equation}
We also have
parameter perturbations of the form
\begin{equation}
J \in [J_0/\gamma_p,\gamma_pJ_0],\, K \in [K_0/\gamma_p,\gamma_pK_0],
\label{eq:motor_position_params}
\end{equation}
with $\gamma_p = 4$.

Using our algorithm we find that we simultaneously stabilise these
under a common Lyapunov function
using output feedback law
\begin{equation}
 u = 
\begin{bmatrix}
    -19.7 & -18.7
\end{bmatrix}
y
\label{eq:motor_position_feedback}
\end{equation}

The identified polyhedron is also shown in \cref{fig:motor} 
and has only 9 vertices.

\section{Conclusions and Future Work}
We presented an algorithm for finding
polyhedral Lyapunov functions of fixed complexity
for feedback design.

We believe that this can be used to tackle interesting research directions.
We plan to extend the theory to open systems and interconnections,
and apply it to the nonlinear contraction setting.


\bibliographystyle{abbrv}

\bibliography{21_CDC,example_websites,fulvio}

\begin{thebibliography}{10}

\bibitem{ctms_control_tutorials}
University of {Michigan} control {{Tutorials}} for {{MATLAB}} and {{Simulink}}.
\newblock \url{https://ctms.engin.umich.edu/CTMS/index.php?aux=Home}.

\bibitem{ambrosino_convex_2012}
R.~Ambrosino, M.~Ariola, and F.~Amato.
\newblock A {{Convex Condition}} for {{Robust Stability Analysis}} via
  {{Polyhedral Lyapunov Functions}}.
\newblock 50(1):490--506.

\bibitem{blanchini_set-theoretic_2015}
F.~Blanchini and S.~Miani.
\newblock {\em Set-Theoretic Methods in Control}.
\newblock Systems \& Control: Foundations \& Applications. {Birkhäuser}, 2
  edition.

\bibitem{boyd_linear_1994}
S.~Boyd, L.~El~Ghaoui, E.~Feron, and V.~Balakrishnan.
\newblock {\em Linear {{Matrix Inequalities}} in {{System}} and {{Control
  Theory}}}.
\newblock Studies in {{Applied}} and {{Numerical Mathematics}}. {Society for
  Industrial and Applied Mathematics}.

\bibitem{boyd_convex_2004}
S.~Boyd and L.~Vandenberghe.
\newblock {\em Convex {{Optimization}}}.
\newblock {Cambridge University Press}.

\bibitem{queiroz_lyapunov_2004}
F.~Clarke.
\newblock Lyapunov {{Functions}} and {{Feedback}} in {{Nonlinear Control}}.
\newblock In M.~S. Queiroz, M.~Malisoff, and P.~Wolenski, editors, {\em Optimal
  {{Control}}, {{Stabilization}} and {{Nonsmooth Analysis}}}, volume 301 of
  {\em Lecture {{Notes}} in {{Control}} and {{Information Sciences}}}, pages
  267--282. {Springer Berlin Heidelberg}.

\bibitem{nash_origins_1990}
G.~B. Dantzig.
\newblock Origins of the simplex method.
\newblock In S.~G. Nash, editor, {\em A History of Scientific Computing}, pages
  141--151. {ACM}.

\bibitem{dontchev_implicit_2009}
A.~L. Dontchev and R.~T. Rockafellar.
\newblock {\em Implicit {{Functions}} and {{Solution Mappings}}: {{A View}}
  from {{Variational Analysis}}}.
\newblock Springer {{Monographs}} in {{Mathematics}}. {Springer New York}.

\bibitem{forni_differential_2014}
F.~Forni and R.~Sepulchre.
\newblock A {{Differential Lyapunov Framework}} for {{Contraction Analysis}}.
\newblock 59(3):614--628.

\bibitem{guglielmi_polytope_2017-1}
N.~Guglielmi, L.~Laglia, and V.~Protasov.
\newblock Polytope {{Lyapunov Functions}} for {{Stable}} and for {{Stabilizable
  LSS}}.
\newblock 17(2):567--623.

\bibitem{kousoulidis_optimization_2020}
D.~Kousoulidis and F.~Forni.
\newblock An {{Optimization Approach}} to {{Verifying}} and {{Synthesizing
  K}}-cooperative {{Systems}}.
\newblock In {\em 21st {{IFAC World Congress}} ({{IFAC}}-{{V}} 2020)}.

\bibitem{lohmiller_contraction_1998-1}
W.~Lohmiller and J.~Slotine.
\newblock On {{Contraction Analysis}} for {{Non}}-linear {{Systems}}.
\newblock 34(6):683--696.

\bibitem{mcmullen_maximum_1970}
P.~McMullen.
\newblock The maximum numbers of faces of a convex polytope.
\newblock 17(2):179--184.

\bibitem{miani_maxis-g:_2005}
S.~Miani and C.~Savorgnan.
\newblock {{MAXIS}}-{{G}}: A software package for computing polyhedral
  invariant sets for constrained {{LPV}} systems.
\newblock In {\em Proceedings of the 44th {{IEEE Conference}} on {{Decision}}
  and {{Control}}}, pages 7609--7614.

\bibitem{nocedal_numerical_2006}
J.~Nocedal and S.~J. Wright.
\newblock {\em Numerical Optimization}.
\newblock Springer Series in Operations Research. {Springer}, 2nd ed edition.

\bibitem{Pavlov2005}
A.~Pavlov, N.~van~de Wouw, and H.~Nijmeijer.
\newblock {\em Uniform Output Regulation of Nonlinear Systems: A Convergent
  Dynamics Approach}.
\newblock Systems {\&} Control: Foundations {\&} Applications. Birkhäuser,
  2005.

\bibitem{polanski_absolute_2000-1}
A.~Polański.
\newblock On absolute stability analysis by polyhedral {{Lyapunov}} functions.
\newblock 36(4):573--578.

\bibitem{Russo2010}
G.~Russo, M.~Di~Bernardo, and E.~Sontag.
\newblock Global entrainment of transcriptional systems to periodic inputs.
\newblock {\em PLoS Computational Biology}, 6(4):e1000739, 04 2010.

\end{thebibliography}

\end{document}